\documentclass[a4paper]{amsart}
\usepackage[latin1]{inputenc}
\usepackage{amssymb}
\usepackage{amsmath}
\usepackage{mathrsfs}
\usepackage{eufrak}
\usepackage{amsthm}
\usepackage{amsfonts}
\usepackage{textcomp}
\usepackage{graphicx}
\usepackage[pdftex]{color}
\usepackage{paralist}
\usepackage[shortlabels]{enumitem}
\usepackage[hyperpageref]{backref}
\usepackage[pagebackref]{hyperref}
\renewcommand*{\backref}[1]{}  
   \renewcommand*{\backrefalt}[4]{
      \ifcase #1 
         Not cited.
      \or
         Cited on page #2.
      \else
         Cited on pages #2.
      \fi}
\usepackage{comment}
\usepackage[arrow, matrix, curve]{xy}

\usepackage{adjustbox}

\usepackage{tikz} \usetikzlibrary{matrix,patterns,shapes,decorations.pathmorphing,decorations.pathreplacing,calc,arrows} 
\usepackage{tikz-cd} 
\usepackage{tkz-graph}

\newtheorem*{corollary*}{Corollary}
\newtheorem{theorem}{Theorem}[section]

\newtheorem{question}[theorem]{Question}

\newtheorem*{claim*}{Claim}

\theoremstyle{definition}
\newtheorem{definition}[theorem]{Definition}

\newtheorem{main conjecture}[theorem]{Main Conjecture}
\newtheorem*{theorem }{Theorem}

\newtheorem{problem}[theorem]{Problem}
\newtheorem{example}[theorem]{Example}

\theoremstyle{remark}

\numberwithin{equation}{section}

\makeatletter
\renewcommand*\env@matrix[1][\
arraystretch]{%
  \edef\arraystretch{#1}%
  \hskip -\arraycolsep
  \let\@ifnextchar\new@ifnextchar
  \array{*\c@MaxMatrixCols c}}
\makeatother

\renewcommand{\mod}{\operatorname{mod}}

\newcommand{\Ext}{\operatorname{Ext}}
\newcommand{\End}{\operatorname{End}}
\newcommand{\CM}{\operatorname{CM}}

\newcommand{\fpdim}{\operatorname{fpdim}}

\newcommand{\Hom}{\operatorname{Hom}}
\newcommand{\add}{\operatorname{\mathrm{add}}}
\renewcommand{\top}{\operatorname{\mathrm{top}}}
\newcommand{\rad}{\operatorname{\mathrm{rad}}}
\newcommand{\soc}{\operatorname{\mathrm{soc}}}

\newcommand{\op}{\operatorname{\mathrm{op}}}

\renewcommand{\mod}{\operatorname{mod}}

\newcommand{\Z}{\mathbb{Z}}

\newcommand{\RHom}{\operatorname{\mathsf{R}Hom}}

\newcommand{\domdim}{\operatorname{domdim}}

\newcommand{\idim}{\operatorname{idim}}

\newcommand{\pdim}{\operatorname{pdim}}
\newcommand{\gldim}{\operatorname{gldim}}
\newcommand{\DynA}{\mathbb{A}}
\newcommand{\DynD}{\mathbb{D}}
\newcommand{\DynE}{\mathbb{E}}

\begin{document}

\title[{Preprojective algebras and generalisations: A short survey}]{Preprojective algebras and generalisations: A short survey}
\date{\today}
\dedicatory{Dedicated to the memory of Idun Reiten}

\subjclass[2010]{16G10, 16E10}

\keywords{Cohen-Macaulay algebras, preprojective algebras, total preprojective algebras, contracted preprojective algebras, homological dimensions}

\author[Chan]{Aaron Chan}
\address[Chan]{Graduate School of Mathematics, Nagoya University, Furocho, Chikusaku, Nagoya 464-8602, Japan}
\email{aaron.kychan@gmail.com}

\author[Iyama]{Osamu Iyama}
\address[Iyama]{Graduate School of Mathematical Sciences, University of Tokyo, 3-8-1 Komaba Meguro-ku Tokyo 153-8914, Japan}
\email{iyama@ms.u-tokyo.ac.jp}

\author[Marczinzik]{Ren\'{e} Marczinzik}
\address[Marczinzik]{Mathematical Institute of the University of Bonn, Endenicher Allee 60, 53115 Bonn, Germany}
\email{marczire@math.uni-bonn.de}

\begin{abstract}
The preprojective algebra of a hereditary algebra $H$ can be defined as a certain orbit construction of the regular representation generated by the Auslander-Reiten translation.  In this short survey, we will look at two important generalisations, namely, the contracted preprojective algebra and the total preprojective algebra. We will include several open problems and questions motivated by examples in the hope to stimulate future research on general orbit algebras of $H$-modules.
\end{abstract}

\maketitle
 
 \section*{Introduction}
Preprojective algebras have become classical in representation theory, with deep ties to quiver varieties \cite{Nak94}, quantum groups \cite{KS97, Lus91}, cluster algebras \cite{GLS2}, Kleinian singularities \cite{Cra00,CBH}, and Cohen-Macaulay modules \cite{Aus1, GL91,Sko}.

Their most common definition due to Gelfand and Ponomarev \cite{GP} goes as follows.
Let $Q=(Q_0,Q_1)$ be a finite quiver with vertex set $Q_0$ and arrow set $Q_1$. Then define the \emph{double quiver}
$\overline{Q}$ as the quiver obtained from $Q$ by adding an arrow $\overline{\alpha}: j \rightarrow i$ for every arrow $\alpha:i \rightarrow j$ in $Q$.
The \emph{preprojective algebra} associated to $Q$ is the algebra 
\[\Pi(Q)=K \overline{Q}/\langle \rho \rangle\ \mbox{ with }\ \rho:= \sum\limits_{\alpha \in Q_1}^{}{(\alpha \overline{\alpha}- \overline{\alpha} \alpha)},\]
which does not depend on the orientation of the quiver $Q$.
We survey some of the most important results on preprojective algebras, as well
as connections to other fields, mainly concerning representation theory of finite-dimensional
algebras. Then we will focus on the following generalisation of preprojective algebras with recent applications, which is the main object of study in this survey.

\begin{definition}\cite{CIM2} \label{maindefinition}
Let $K$ be a field, and let $H$ be a finite dimensional hereditary $K$-algebra. For a finitely generated $H$-module $X$, the \emph{$X$-preprojective algebra} is defined as an orbit algebra
\[\Psi_H^X:=\bigoplus\limits_{i \geq 0}(\Psi_H^X)_i\ \mbox{ with }\ (\Psi_H^X)_i := \Hom_H(X, \tau^{-i}(X)).\]
\end{definition}

Here $\tau^{-1}$ denotes the inverse  \emph{Auslander-Reiten translation} given by $\tau^{-1} =\Ext_H^1(DH,-) \colon \mod H \rightarrow \mod H$ with $D=\Hom_K(-,K)$ the natural duality. The multiplication is defined by $g\cdot f:=\tau^{-i}(g)f \in(\Psi_H^X)_{i+j}$ for $f\in(\Psi_H^X)_i$ and $g\in (\Psi_H^X)_j$.
Here we can always assume that, up to Morita equivalence, $X$ is a basic module, see \cite[Proposition 2.1]{CIM2}.
It is well known that, by a classical result of Ringel \cite{R}, the classical preprojective algebra $\Pi(Q)$, defined by quiver-and-relation above, coincides with the case when $H=KQ$ and $X=H$. 

It seems that $\Psi_H^X$ has not been studied systematically apart from the case when $X=H$. 
In this note,
we survey some nice properties for the following two other choices of $X$ when $H$ is a path algebra of Dynkin type.
\begin{enumerate}
    \item $X$ is projective. In this case, $\Psi_H^X$ is Morita equivalent to the subalgebra $e \Pi(H) e$ for an idempotent $e$ of $\Pi(H)$, and we call $\Psi_H^X$ a \emph{contracted preprojective algebra}.
    \item $X$ is the direct sum of all indecomposable $H$-modules. In this case, we call $\Psi_H^X$ the \emph{total preprojective algebra} of $H$.
\end{enumerate}
To the best of our knowledge, the first detailed analysis of contracted preprojective algebras appears in the work of Iyama and Wemyss in \cite{IW}. Recently, we used contracted preprojective algebras in \cite{CIM} to give the first non-trivial classes of Cohen-Macaulay algebras (in the sense of Auslander and Reiten; see Definition \ref{def:CM alg}), and also to give an answer to an open question (Question \ref{Auslander-Reiten question}) due to Auslander and Reiten \cite{AR2}.

Total preprojective algebras first appeared in the work of Geiss, Leclerc and Schr\"oer 
as the endomorphism algebras of certain 2-cluster tilting modules for preprojective algebras of Dynkin type $Q$, see \cite{GLS1,GLS2}. 
It was shown in \cite{CIM2} that these algebras are canonically isomorphic to our total preprojective algebras (Theorem \ref{main kQ}).

The goal of this survey article is to present our recent results on contracted preprojective algebras and total preprojective algebras, and to pose some related problems and questions with motivating examples.
Section 1 provides preliminaries and a brief survey of some of the most important results on preprojective algebras of Dynkin type that are related to the classical representation theory of finite-dimensional algebras.
Section 2 deals with Cohen-Macaulay algebras in the sense of Auslander and Reiten and with contracted preprojective algebras of Dynkin type.
In Section 3, we look at total preprojective algebras of Dynkin type and give an explicit description via quivers-with-relations, using the Auslander-Reiten quiver of the module category.
In Section 4, we discuss the generalisation of total preprojective algebras to higher dimensions, in analogy with the relationship between higher Auslander algebras and classical Auslander algebras.
While various results generalise smoothly when transferring to the higher-dimensional setting for total preprojective algebras, it remains mysterious and complicated to us how higher-dimensional versions of contracted preprojective algebras behave.
In the final Section 5, we pose several questions and problems, hoping that they will motivate further research on the various generalised versions of preprojective algebras explored in this survey.

\section{Background and important results on classical preprojective algebras}
Unless otherwise stated, we deal with finite dimensional $K$-algebras with $K$ a field, and assume that modules are finitely generated right modules.
We will assume that $K$-algebras are basic unless stated otherwise. Note that we can assume this up to Morita equivalence.
For background on representation theory and homological algebra of finite dimensional algebras we refer, for example, to the textbooks \cite{ARS,ASS,Kra,SkoYam}. 
We follow the conventions as in \cite{ASS,SkoYam}
on multiplication in path algebras, i.e. the path $ab$ means first $a$, then $b$. 
Before giving some theoretical results on preprojective algebras we start with a simple example:
\begin{example}\label{eg:Pi(D5)}
Let $Q$ be the Dynkin quiver of type $D_5$ with the following orientation:
\[\begin{tikzcd}
	& 2 \\
	1 & 3 & 4 & 5
	\arrow["a", from=1-2, to=2-2]
	\arrow["b"', from=2-1, to=2-2]
	\arrow["c"', from=2-2, to=2-3]
	\arrow["d"', from=2-3, to=2-4]
\end{tikzcd}\]
Then the preprojective algebra 
$\Pi(Q)=K\overline{Q}/I$ is given by the quiver
\[ \overline{Q}= \quad \begin{tikzcd}
	& 2 \\
	1 & 3 & 4 & 5
	\arrow["a", shift left, from=1-2, to=2-2]
	\arrow["b"', shift right, from=2-1, to=2-2]
	\arrow["{\overline{a}}", shift left, from=2-2, to=1-2]
	\arrow["{\overline{b}}"', shift right, from=2-2, to=2-1]
	\arrow["c"', shift right, from=2-2, to=2-3]
	\arrow["{\overline{c}}"', shift right, from=2-3, to=2-2]
	\arrow["d"', shift right, from=2-3, to=2-4]
	\arrow["{\overline{d}}"', shift right, from=2-4, to=2-3]
\end{tikzcd}\]
with 
relations $I=\langle a \overline{a}, b \overline{b}, -\overline{a} a-\overline{b} b+c\overline{c}, \overline{c} c-d \overline{d}, \overline{d} d  \rangle.$
\end{example}

\subsection{Basic properties and description}

We collect now the most important results on classical preprojective algebras from the point of view of representation theory of finite dimensional algebras.
We start with the classical characterisations of preprojective algebras for hereditary algebras $H=KQ$ with $Q$ a finite acyclic quiver.
\begin{theorem}[\cite{R}, see also \cite{Sa}]\label{ppa 3 defs}
Let $H=KQ$ be a finite dimensional hereditary $K$-algebra. Then the following algebras are isomorphic:
\begin{enumerate}[\rm(1)]
    \item The preprojective algebra $\Pi(Q)$ of $Q$.
    \item The $H$-preprojective algebra $\Psi_H^H=\bigoplus\limits_{i \geq 0}^{}{\Hom_H (H, \tau^{-i}(H)})$.
    \item The tensor algebra of the $H$-bimodule $\Ext_H^1(DH,H)$ over $H$.
\end{enumerate}
\end{theorem}
In view of Theorem \ref{ppa 3 defs}, for any hereditary $K$-algebra $H$ (which is not necessarily the path algebra of a quiver), we define the \emph{preprojective algebra} of $H$ by $\Pi(H):=\Psi_H^H$.
The following conditions are equivalent (see \cite{ARS}):
\begin{enumerate}[\rm(i)]
\item $H$ is representation-finite.
\item The underlying valued graph of $H$ is Dynkin, i.e. one of $\DynA_n, \mathbb{B}_n, \mathbb{C}_n,\DynD_n, \DynE_6, \DynE_7, \DynE_8, \mathbb{F}_4, \mathbb{G}_2$.
\item $\Pi(H)$ is finite dimensional.
\end{enumerate}
For connected acyclic quivers over an algebraically closed field, this fits into the
usual finite--tame--wild trichotomy. In Dynkin type, the path algebra $KQ$ is
representation-finite and the preprojective algebra $\Pi(Q)$ is finite dimensional.
In extended Dynkin type, $KQ$ is tame and $\Pi(Q)$ is infinite dimensional but
still highly controlled; more generally, deformed preprojective algebras of extended
Dynkin quivers are prime noetherian rings of Gelfand--Kirillov dimension two
\cite[Theorem 0.8]{CB2}. This is the case closely related to Kleinian singularities
via the Crawley-Boevey--Holland construction \cite{CBH}. For all other connected
acyclic quivers, $KQ$ is wild; then $\Pi(Q)$ is again infinite dimensional, and for
connected non-Dynkin $Q$ it is $2$-Calabi--Yau in the bimodule sense
\cite[Theorem 1.1]{CK}.
In this survey, we will focus mostly in the case when $H=KQ$ for $Q$ of (simply-laced) Dynkin type.

We now list some nice structural properties of preprojective algebras of Dynkin type.
Recall that a finite dimensional $K$-algebra $A$ is called \emph{selfinjective algebra} if the right $A$-module $A_A$ is injective. We refer to the textbook \cite{SkoYam} for an in-depth introduction to selfinjective algebras.  Suppose that $P_1, P_2,...,P_n$ are the indecomposable projective $A$-modules of a selfinjective algebra $A$, then the \emph{Nakayama permutation} of $A$ is the permutation $\pi: \{1,...,n\} \rightarrow \{1,...,n\}$ defined by $\top P_i \cong \soc P_{\pi(i)}$.

For each Dynkin graph $\Delta$, denote by $h_\Delta$ the Coxeter number of type $\Delta$ (see table below).
There is a natural involution $\nu_\Delta$ on $\Delta$, which is non-trivial only in the $\DynA\DynD\DynE$ cases, which can be described by swapping appropriate vertices of $\Delta$.
To describe them explicitly we fix the following enumeration of vertices of Dynkin graphs:
\begin{align}
\DynA_n: & \xymatrix{ 1 \ar@{-}[r] & 2 \ar@{-}[r] & \cdots \ar@{-}[r]& n-1 \ar@{-}[r] & n} \notag \\
\DynD_n: & \xymatrix@R=10pt{ 1 \ar@{-}[r] & 2 \ar@{-}[r] & \cdots \ar@{-}[r]& n-2\ar@{-}[r]\ar@{-}[d] & n-1 \\  & & &n& } \notag \\
\DynE_6: & \xymatrix@R=10pt{ 1 \ar@{-}[r] & 2 \ar@{-}[r] & 3 \ar@{-}[d]\ar@{-}[r] & 4\ar@{-}[r] & 5 & \\  & & 6&&& } \notag 
\end{align}
Then the Coxeter number $h_\Delta$ and the involution $\nu_\Delta$ can be summarised as follows
\[
\begin{array}{c||c|c|c|c|c|c|c|c}
\Delta     & \DynA_n & \mathbb{B}_n,\mathbb{C}_n & \DynD_n & \DynE_6 & \DynE_7 & \DynE_8 & \mathbb{F}_4 & \mathbb{G}_2 \\ \hline
h_\Delta   & n+1     & 2n                        & 2n-2    & 12      & 18      & 30      & 12 & 6 \\
\nu_\Delta & i \leftrightarrow n+1-i & \text{identity}&  n-1 \leftrightarrow n \text{, if $n$ odd} & 1\leftrightarrow 5 & \text{identity} & \text{identity}& \text{identity}& \text{identity}\\
 &  & &  \text{identity, if $n$ even} & 2\leftrightarrow 4 & & & & 
\end{array}
\]

We can now state important properties of preprojective algebras of Dynkin type:
\begin{theorem}
Let $H$ be a hereditary $K$-algebra of Dynkin type $\Delta$.
Then $\Pi(H)$ is a selfinjective algebra with Nakayama permutation $\nu_\Delta$.
The Loewy length of $\Pi(H)$ is given by $h_\Delta-1$.
\end{theorem}
We refer to \cite{G,Hu} for more details and a modern proof.

Recall that an $A$-module $M$ is called \emph{periodic} of period $n$ if $\Omega_A^n(M) \cong M$ and $n>0$ is minimal with this property.
The algebra $A$ is called \emph{periodic} if $A$ as an $A$-bimodule is periodic. In this case, it is known that $A$ is a selfinjective algebra and $\Omega_A^n(M) \cong M$ in the stable module category $\underline{\mod} A$ for every $A$-module $M$, where $n$ is the period of $A$ as an $A$-bimodule. We refer to \cite[Chapter IV]{SkoYam} for the definition of the stable module category and more information on periodic modules and algebras. 
It turns out that preprojective algebras provide examples of periodic algebras.
\begin{theorem}
Let $H=KQ$ be a path algebra of $\DynA\DynD\DynE$ type. Then the preprojective algebra $\Pi(H)$ is periodic of period dividing 6.
\end{theorem}
The first published proof of this theorem seems to be due to Erdmann and Snashall \cite{ErSn}, but it is also well-known that this was first shown by Ringel and Schofield in an unpublished note.

\subsection{Cluster tilting theory}

Closely related to the periodicity, preprojective algebras of Dynkin type exhibit another important property -- being stably 2-Calabi-Yau -- which has deep connection to cluster tilting theory; see \cite[7.5]{GLS1}.
Recall that a $K$-linear Hom-finite triangulated category $\mathcal{C}$ with shift functor $\Sigma$ is called \emph{2-Calabi-Yau} if for all objects $X$ and $Y$ of $\mathcal{C}$ there is a bifunctorial isomorphism $\Hom_{\mathcal{C}}(X, \Sigma Y) \cong D \Hom_{\mathcal{C}}(Y, \Sigma X)$, where $D=\Hom_K(-,K)$ is the natural vector space duality.  Note that the stable module category of any selfinjective finite dimensional algebra is a triangulated category with shift functor given by $\Omega^{-1}$, we refer for example to \cite[Chapter 6]{Kra}.
\begin{theorem}\cite[7.5]{GLS1}
Let $H=KQ$, where $Q$ is a Dynkin quiver of type $\Delta$. Then the stable
module category of $\Pi(H)$ is $2$-Calabi--Yau.
\end{theorem}

As mentioned, the 2-Calabi-Yau property is deeply related to cluster tilting theory, which was motivated by the study of cluster algebras; this connection is first shown by Geiss, Leclerc and Schr\"oer in \cite{GLS2}.  Let us recall this important structural result for preprojective algebras of Dynkin type.
Recall that an $A$-module $M$ is called \emph{rigid} if $\Ext_A^1(M,M)=0$ and \emph{maximal rigid} if $M$ is rigid and $M \oplus N$ is not rigid for any indecomposable $A$-module $N$ that is not a direct summand of $M$. 
The following is a central result in \cite{GLS2}:

\begin{theorem}\cite[Theorem 2.1]{GLS2}\label{num summand of rigid}
Let $H=KQ$, where $Q$ is a Dynkin quiver of type $\Delta$, and let $M$ be a
rigid $\Pi(H)$-module. Then the number of indecomposable non-isomorphic direct
summands of $M$ is bounded by the number of positive roots of $\Delta$.
\end{theorem}

This follows from a deeper result on 2-cluster tilting modules over $\Pi(KQ)$.
Recall that an $A$-module $M$ is called \emph{$d$-cluster tilting} for $d \geq 1$ if 
\begin{eqnarray*}
    \add M&=& \{ Y \in \mod A \mid \Ext_A^i(Y,M)=0 \ \text{for} \ i=1,...,d-1 \}\\
    &=& \{ Y \in \mod A \mid \Ext_A^i(M,Y)=0 \ \text{for} \ i=1,...,d-1 \}.
\end{eqnarray*}
This notion was first introduced by Iyama in \cite{I} and developed into its own subarea of modern homological algebra and representation theory, called higher Auslander-Reiten theory. 
In \cite{A}, Auslander established what is now called Auslander correspondence, which is a correspondence between representation-finite algebras and the so-called Auslander algebras $A$, i.e. those that satisfy $\gldim A \leq 2 \leq \domdim A$, see also \cite[Chapter VI.5]{ARS}.
The Auslander correspondence was generalised by Iyama in \cite{I} to the higher Auslander correspondence. Recall that an algebra $A$ is called \emph{$d$-Auslander algebra} for some $d \geq 1$ if $\gldim A \leq d+1 \leq \domdim A.$
The higher Auslander correspondence then states that there is a bijective correspondence between $d$-Auslander algebras $A$ and $d$-cluster tilting $B$-modules $M$ of an algebra $B$, where $A$ is given by $A=\End_B(M)$, we refer for example to \cite{I} for more details on this correspondence and a proof.

Each direct summand of a 2-cluster tilting module is rigid; in fact, an $A$-module $M$ is maximal rigid if and only if $M$ is a 2-cluster tilting $A$-module in the setting of preprojective algebras of Dynkin type.  Now Theorem \ref{num summand of rigid} can be seen as an immediate consequence of the following result.

\begin{theorem}\cite[Theorem 2.2]{GLS2} \cite[Theorem 1]{GLS1}\label{GLS 2CT over Pi}
Let $H=KQ$, where $Q$ is a Dynkin quiver of type $\Delta$. Then $\Pi(H)$ has a
$2$-cluster tilting module, and the number of non-isomorphic indecomposable
summands of each $2$-cluster tilting module is given by the number of positive
roots of $\Delta$.
\end{theorem}

We give a simple example:
\begin{example}
Let $A=KQ/I=\Pi(\DynA_3)$ be the preprojective algebra of Dynkin type $\DynA_3$, i.e.
\[ Q = \begin{tikzcd}
	1 & 2 & 3
	\arrow["a"', shift right, from=1-1, to=1-2]
	\arrow["{\overline{a}}"', shift right, from=1-2, to=1-1]
	\arrow["b"', shift right, from=1-2, to=1-3]
	\arrow["{\overline{b}}"', shift right, from=1-3, to=1-2]
\end{tikzcd}\;\; \text{ and } \;\; I=\langle a\overline{a},\overline{a}a-b\overline{b}, \overline{b}b \rangle.\]
$A$ is representation-finite with 12 indecomposable modules. As there are 6 positive roots, every 2-cluster tilting $A$-module has 6 indecomposable non-isomorphic direct summands.
There are 14 2-cluster tilting modules up to isomorphism, we leave it to the interested reader to classify all of them.  In fact, the stable module category of this algebra is equivalent to the cluster category of type $\DynA_3$, which explains the number 14, see for example  \cite[Section 1.2, 1.3]{A} for more details on this example. We give one specific example of a 2-cluster tilting module:
Let $S_i$ denote the simple $A$-modules and $P_i$ the indecomposable projective $A$-modules.
Then $M:=A \oplus S_1 \oplus S_3 \oplus \rad(P_2)$ is a 2-cluster tilting module.
The endomorphism ring of $M$, which is a higher Auslander algebra of global dimension 3, is given by the following quiver 
\[
Q'=\begin{tikzcd}
	S_3 & \rad(P_2) & S_1 \\
	P_1 & P_2 & P_3
	\arrow["{a_5}"', from=1-1, to=2-1]
	\arrow["{a_3}"', from=1-2, to=1-1]
	\arrow["{a_2}", from=1-2, to=1-3]
	\arrow["{a_4}", from=1-2, to=2-2]
	\arrow["{a_1}", from=1-3, to=2-3]
	\arrow["{a_6}"{description}, from=2-1, to=1-2]
	\arrow["{a_7}", from=2-2, to=2-1]
	\arrow["{a_8}"', from=2-2, to=2-3]
	\arrow["{a_9}"{description}, from=2-3, to=1-2]
\end{tikzcd}\]
with relations
$I'=\langle a_1 a_9, a_3 a_5-a_4 a_7, a_2 a_1-a_4 a_8, a_5 a_6, a_6 a_3, a_7 a_6-a_8 a_9, a_9 a_2 \rangle.$
\end{example}

The previous Theorem \ref{GLS 2CT over Pi} opened a new connection to the theory of cluster algebras, we refer for example to the survey articles \cite{GLS} and \cite{GLS3} for more information.
It turns out that there is a canonical choice of 2-cluster tilting module for $\Pi(KQ)$, whose existence was first noticed in \cite{GLS1}.  This will come up again in the context of total preprojective algebras in the upcoming Section 3.

\subsection{Connection to Cohen-Macaulay representations}

There is a fourth important description of preprojective algebras coming from singularity theory and commutative algebra.
Let $k$ be an algebraically closed field of characteristic $0$. A \emph{simple singularity} over $k$ is a hypersurface of the form
\[
R \;=\; k[[x_0,x_1,\dots,x_d]]/(f_\Delta^{\,d}),
\]
where $\Delta\in\{\DynA_n,\DynD_n,\DynE_6,\DynE_7,\DynE_8\}$ is a Dynkin diagram and $f_\Delta^{\,d}$ is one of:
\begin{align*}
f_{\DynA_n}^{\,d} &= x_0^2 + x_1^{n+1} + x_2^2 + \cdots + x_d^2, & n\ge 1,\\
f_{\DynD_n}^{\,d} &= x_0^2 x_1 + x_1^{n-1} + x_2^2 + \cdots + x_d^2, & n\ge 4,\\
f_{\DynE_6}^{\,d} &= x_0^3 + x_1^{4} + x_2^2 + \cdots + x_d^2,\\
f_{\DynE_7}^{\,d} &= x_0^3 + x_0 x_1^{3} + x_2^2 + \cdots + x_d^2,\\
f_{\DynE_8}^{\,d} &= x_0^3 + x_1^{5} + x_2^2 + \cdots + x_d^2.
\end{align*}
Simple singularities are exactly the hypersurface singularities of finite deformation type; see \cite{Ar1,Ar2}.

Recall that a complete local Cohen-Macaulay ring $R$ is said to have \emph{finite Cohen-Macaulay type} if there are only finitely many indecomposable maximal Cohen-Macaulay $R$-modules up to isomorphism. The classification of Cohen-Macaulay rings of finite Cohen-Macaulay type is wide open, see for example the survey \cite{Y}. The theorem of Buchweitz-Greuel-Schreyer and Kn\"orrer \cite{BGS,K} provides a representation-theoretic characterization of simple singularities in the Gorenstein case:

\begin{theorem}
Let $R$ be a complete local Gorenstein ring of Krull dimension $d$ over an algebraically closed field of characteristic $0$. Then $R$ has finite Cohen-Macaulay type if and only if $R$ is a simple singularity.
\end{theorem}

See also the textbooks \cite{Y,LW} for further details.
The connection of simple singularities with preprojective algebras of Dynkin type can be highlighted as follows.
\begin{theorem} \label{simple sing vs ppa}
Let $d \geq 2$ be even and $R=k[[x_0,x_1,...,x_d]]/(f_\Delta^{d})$ be a simple singularity over an algebraically closed field $k$ of characteristic 0. Let $M$ be the direct sum of all indecomposable maximal Cohen-Macaulay $R$-modules. Then the stable endomorphism ring $\underline{\End}_R(M)$ is isomorphic to the preprojective algebra of the corresponding Dynkin type $\Delta$.
\end{theorem}

For more details and explanations, we refer to the last section of the survey \cite{Sko} and the article \cite{Sol}.
In case of odd $d$, the so-called twisted preprojective algebras appear as stable endomorphism rings, we refer to the last section of \cite{Sko} for more details on the odd $d$ case.

Notice that the ordinary endomorphism ring $\End_R(M)$ is isomorphic to the completion of the preprojective algebra of extended Dynkin type $\widetilde{\Delta}$, see \cite[section 5.1]{I4} and references therein.

\subsection{Outlook on Calabi-Yau completions}
The last description of preprojective algebras is via Calabi-Yau completions, we refer for example to \cite{Kel2} for unexplained notions.
Let $A$ be a dg $K$-algebra, and $A^e:=A\otimes_KA^{\op}$ its enveloping dg algebra.
We call $A$ \emph{smooth} if $A$ belongs to the perfect derived category of $A^e$. For $d\in\Z$, we call $A$ \emph{$d$-Calabi-Yau} if $A$ is smooth and $\RHom_{A^e}(A,A^e)\simeq A[-d]$ in $D(A^e)$.
Let $\Theta_d$ be a projective resolution of $\RHom_{A^e}(A,A^e)[d]$. The \emph{$(d+1)$-derived preprojective algebra} (also called the \emph{$(d+1)$-Calabi-Yau completion}) of $A$ is defined as the tensor dg algebra
\[{\mathbf\Pi}_{d+1}(A):=T_A(\Theta_d).\]
We have the following important result.

\begin{theorem}\cite{Kel}
   For each smooth dg $K$-algebra $A$ and $d\in\Z$, ${\mathbf\Pi}_{d+1}(A)$ is a $(d+1)$-Calabi-Yau dg algebra.
\end{theorem}

Let $Q$ be a connected acyclic quiver and $KQ$ its path algebra. Then we have
\[\Pi(KQ)=H^0({\mathbf\Pi}_2(KQ)).\]
Moreover, if $Q$ is non-Dynkin, then ${\mathbf\Pi}_2(KQ)$ is quasi-isomorphic to $\Pi(KQ)$. This gives an explanation why $\Pi(KQ)$ for non-Dynkin quivers $Q$ is 2-Calabi-Yau in the bimodule sense.
Similar results hold for ${\mathbf\Pi}_{d+1}(A)$ of a finite dimensional algebra $A$ with global dimension at most $d$, see Section 4 and \cite{HIO}.

\subsection{Remarks}

There are many topics related to preprojective algebras that we do not touch on in this short survey article, see also the references in the beginning of the introduction. We just mention three important aspects here at the end for the interested readers.
There are deformed preprojective algebras that appear in various areas such as Kleinian singularities, the Deligne-Simpson problem, integrable systems and noncommutative geometry, see for example \cite{CK}. They also appear in the proof of Kac's theorem that is a generalisation of Gabriel's theorem. For a textbook introduction, we refer for example to \cite[Chapter 8]{DW}.
Another very important aspect of preprojective algebras is the connection to 
cluster algebras and quantum groups via Lusztig's nilpotent varieties that can be interpreted as varieties of modules over the preprojective algebra, we refer to \cite{Ki,KS97,Lus91,GLS1,GLS2,GLS3,GLS4} for more details and connections to topics such as cluster algebras and canonical bases.
Finally, we mention that representation theoretic and homological properties of preprojective algebras (not only of Dynkin type) have nice interactions with combinatorial structures such as lattice theory, we refer for example to \cite{IRRT,ORT,BIRS}.

\section{Contracted preprojective algebras and Cohen-Macaulay algebras in the sense of Auslander and Reiten}
We start with a brief overview on Cohen-Macaulay algebras for finite dimensional algebras.
Recall that an Artin algebra $A$ is called \emph{Iwanaga-Gorenstein} if $\idim A_A= \idim _{A}A < \infty$. If $P^{< \infty}$ denotes the full subcategory of $\mod A$ consisting of modules of finite projective dimension and $I^{< \infty}$ denotes the full subcategory of modules of finite injective dimension, then it is well known that $A$ is Iwanaga-Gorenstein if and only if $P^{< \infty}=I^{< \infty}$. Iwanaga-Gorenstein algebras contain several important classes of algebras such as selfinjective algebras, gentle algebras, cluster-tilted algebras \cite{SkoYam,GR,KR}, and many others.
Auslander and Reiten gave the following generalisation of Iwanaga-Gorenstein algebras in \cite{AR1,AR2}:
\begin{definition}\label{def:CM alg}
A finite dimensional algebra $A$ is called \emph{Cohen-Macaulay} if there is an $A$-bimodule $W$ such that the functor $\Hom_A(W,-)$ induces an equivalence between the categories $I^{< \infty}$ and $P^{< \infty}$. In this case $W$ is called a \emph{dualizing $A$-module}.
\end{definition}
When $A$ is Cohen-Macaulay, then the dualizing module $W$ has several nice properties. It is an \emph{Ext-maximal} cotilting module, that is, it is the maximum in the set of cotilting $A$-modules under the natural order given by Ext-vanishing, see \cite{AR2,CIM} for more details and properties. This also shows that it is essentially unique.

The definition of Cohen-Macaulay for finite dimensional algebras is motivated by the classical notion of Cohen-Macaulay commutative local noetherian rings with a dualizing module, we refer to \cite{AR2} for more details and to \cite{BR} for a generalisation to abelian categories.
From now on, when dealing with finite dimensional algebras, by Cohen-Macaulay algebra we always mean the one given above.

Recall that the finitistic projective dimension of a finite dimensional algebra is defined as $\operatorname{fpdim} A= \sup\{ \pdim M \mid M \in \mod A, \ \pdim M< \infty \}$.
An important consequence for an algebra $A$ being Cohen-Macaulay is that its finitistic (left and right) dimension is equal to the injective dimension $\idim W$ of $W$.
Auslander and Reiten introduced this notion of Cohen-Macaulay algebras for finite dimensional algebras in the early 1990's, but there are very few known examples, which are fully listed below.
\begin{enumerate}
    \item Iwanaga-Gorenstein algebras. These are exactly the Cohen-Macaulay algebras with $W=A$.
    \item Algebras $A$ whose left and right finitistic dimension is equal to 0. This is the case if and only if one can choose $W=D A$. Examples include local finite dimensional algebras.
    \item Tensor products of the algebras in case (1) and (2).
\end{enumerate}

A module $M$ over a Cohen-Macaulay algebra $A$ is called \emph{maximal Cohen-Macaulay} if $M \in \CM A:=\{ X \in \mod A \mid \Ext_A^i(X,W)=0 \ \text{for} \ i>0\}$.
This generalises the classical maximal Cohen-Macaulay modules, also known as Gorenstein projective modules, over Iwanaga-Gorenstein algebras.  We refer to \cite[Chapter 6]{Kra} for an introduction to Gorenstein homological algebra.
Auslander and Reiten posed the following question in \cite{AR2}:
\begin{question} \label{Auslander-Reiten question} (Auslander-Reiten)
Let $A$ be a Cohen-Macaulay finite dimensional algebra with $n:=\idim W>0$ and assume that $\Omega^n(\mod A)=\operatorname{CM} A$. Is $A$ then Iwanaga-Gorenstein?
\end{question}

We will give the first non-trivial (in the sense that in general they are not of the form (1), (2) or (3) as in the list above) examples of Cohen-Macaulay finite dimensional algebras using preprojective algebras.
The main result in \cite{CIM} is as follows:
\begin{theorem} \label{CMmaintheorem} \cite[Theorem 4.2]{CIM}
Let $A=\Psi_{KQ}^M$ for $Q$ of $\DynA\DynD\DynE$-type and $M$ a projective $KQ$-module. Then $A$ is a Cohen-Macaulay algebra with $\idim W \leq 2.$
\end{theorem}
Note that, if we let $e$ to be the sum of primitive idempotents corresponding to the indecomposable direct summands of a projective module $M$, then $\Psi_H^M$ is Morita equivalent $e\Pi(H)e$.
Explicit homological dimensions of contracted preprojective algebras such as the global dimension and the dominant dimension can be found in \cite[Theorem 4.6]{CIM} and an explicit description of the dualizing modules can be found in \cite[Proposition 3.11]{CIM}.
We can describe the category of maximal Cohen-Macaulay modules explicitly in important situations due to the next result.
Recall here that an algebra $A$ with minimal injective coresolution 
$0 \rightarrow A \rightarrow I^0 \rightarrow I^1 \rightarrow \cdots$
is called \emph{$d$-Gorenstein} if $\pdim I^i \leq i$ for $i=0,1,..,d-1$, we refer for example to \cite{AR3,FGR,KM} for more information on this class of algebras.
\begin{theorem} \label{syzygyCMtheorem} \cite[Corollary 2.2]{CIM}
Let $A$ be a Cohen-Macaulay algebra with $\idim W=d$. If A is $d$-Gorenstein, then $\operatorname{CM} A=\Omega^d(\mod A)$. In this case we have 
$$\add W =\add\big(\Big(\bigoplus_{0\le i\le d-1}P_i(DA)\Big) \oplus \Omega^d(DA)\big),$$
where $P_i(DA)$ denotes the $i$-th term in a minimal projective resolution of $DA$.
\end{theorem}

We now give a non-trivial example that gives a negative answer to Auslander-Reiten's Question \ref{Auslander-Reiten question}:
\begin{example}
Let $H=KQ$ with $Q$ of Dynkin type $D_5$ as in Example \ref{eg:Pi(D5)}.
Let $e=e_1+e_3+e_4+e_5$ and $A=e\Pi(Q)e$ the contracted preprojective algebra.
Then $A=KQ'/I'$ with 
\[Q'= \begin{tikzcd}
	5 & 4 & 3 & 1
	\arrow["{\overline{d}}", shift left, from=1-1, to=1-2]
	\arrow["{d}", shift left, from=1-2, to=1-1]
	\arrow["{\overline{c}}", shift left, from=1-2, to=1-3]
	\arrow["{c}", shift left, from=1-3, to=1-2]
	\arrow["{\overline{b}}", shift left, from=1-3, to=1-4]
	\arrow["{b}", shift left, from=1-4, to=1-3]
\end{tikzcd}\;\text{ and }\; I'=\langle \overline{d}d, d\overline{d}-\overline{c}c, b\overline{b}, cd\overline{d}\overline{c}-c\overline{c}\overline{b}b- \overline{b}bc\overline{c}\rangle.\]
$A$ is Cohen-Macaulay of $K$-dimension 42 with dualizing module $W=U \oplus \Omega^2(DA)$, where $U$ is the direct sum of all indecomposable projective-injective $A$-modules. We have $d=\idim W=2$ and $A$ is 2-Gorenstein.
By Theorem \ref{syzygyCMtheorem}, we have $\CM A= \Omega^2(\mod A)$. However, $A$ is not Iwanaga-Gorenstein and thus this gives a negative answer to Question \ref{Auslander-Reiten question} of Auslander and Reiten.
\end{example}
More generally, we can fully classify the contracted preprojective algebras $A=e\Pi(KQ)e$ for $Q$ of Dynkin type according to their selfinjective dimension $\idim A$, finitistic dimension, and dominant dimension (\cite[Theorem 4.6]{CIM}).  $A$ satisfies $\CM(A)=\Omega^2(\mod A)$ when both finitistic dimension and dominant dimension are $2$.  In particular, those with infinite selfinjective dimensions and dominant dimension $\geq 2$ all give a negative answer to Question \ref{Auslander-Reiten question}; we refer to \cite[Theorem 4.6]{CIM} for a systematic construction of these algebras.

An important Corollary of Theorem \ref{CMmaintheorem} is the following:
\begin{theorem} \cite[Theorem 4.4]{CIM}
Let $R$ be a simple singularity over an algebraically closed field of characteristic 0 and $M$ a maximal Cohen-Macaulay $R$-module. Then the stable endomorphism ring of $M$ is a finite dimensional Cohen-Macaulay algebra.
\end{theorem}

In general, for a non-projective $KQ$-module $M$, the algebra $\Psi_{KQ}^M$ may no longer be Cohen-Macaulay; see Example \ref{exoticexamples}.
We pose the following problem:
\begin{problem}
Let $H$ be a hereditary algebra of Dynkin type.
For which $M\in \mod H$ is $\Psi_H^M$ Cohen-Macaulay, Iwanaga-Gorenstein or of finite global dimension?
What are other interesting homological or representation-theoretic properties of them?
\end{problem}

In a forthcoming work \cite{CIM3}, we will give a new way to construct non-trivial Cohen-Macaulay algebras by generalising the Iyama-Solberg correspondence \cite{IyaSol} and the (higher) Auslander correspondence.

\section{Total preprojective algebras}
In this section, we look at the total preprojective algebra $\Psi=\Psi_H^\Pi$ of a representation-finite hereditary algebra $H$, where $\Pi$ is the ordinary preprojective algebra of $H$.  Note that here we view $\Pi=\Psi_H^H$ as an $H$-module, which is isomorphic to the direct sum of all indecomposable $H$-modules, one from each isomorphism class.
The starting point for our main theorem in this section is the following result:

\begin{theorem} \cite[Proposition 2.5]{CIM2} \label{explicit2CTobjectpreprojective}
Let $H$ be a representation-finite hereditary algebra with preprojective algebra $\Pi$. Then $\Pi \otimes_{H} \Pi$ is a 2-cluster tilting module in $\mod \Pi$.
\end{theorem}
Thus the endomorphism ring $\End_{\Pi}(\Pi \otimes_{H} \Pi)$ is a higher Auslander algebra and the next theorem tells us that this is exactly the total preprojective algebra of $H$. Here we will denote by $\Pi_1$ the degree one component of the preprojective algebras with respect to the grading as in Definition \ref{maindefinition}.

\begin{theorem}\label{main kQ} \cite[Theorem 1.3]{CIM2}
Let $H$ be a hereditary algebra of Dynkin type with $\Pi$ the preprojective algebra of $H$ and $\Psi$ the total preprojective algebra.
\begin{enumerate}[\rm(a)]
\item We have an isomorphism of algebras
\[\Psi\simeq\End_\Pi(\Pi\otimes_H\Pi).\]
\item $\Psi$ is a 2-Auslander algebra, that is, $\gldim \Psi\le3\le\domdim \Psi$ holds.
\item Let $\Gamma:=\End_H(\Pi)$ be the Auslander algebra of $H$. Then $\Psi$ is isomorphic to the tensor algebra 
\[\Psi\simeq T_\Gamma\Hom_H(\Pi,\Pi\otimes_H\Pi_1).\]
\end{enumerate}
\end{theorem}

We illustrate Theorem \ref{main kQ} with the following figure.
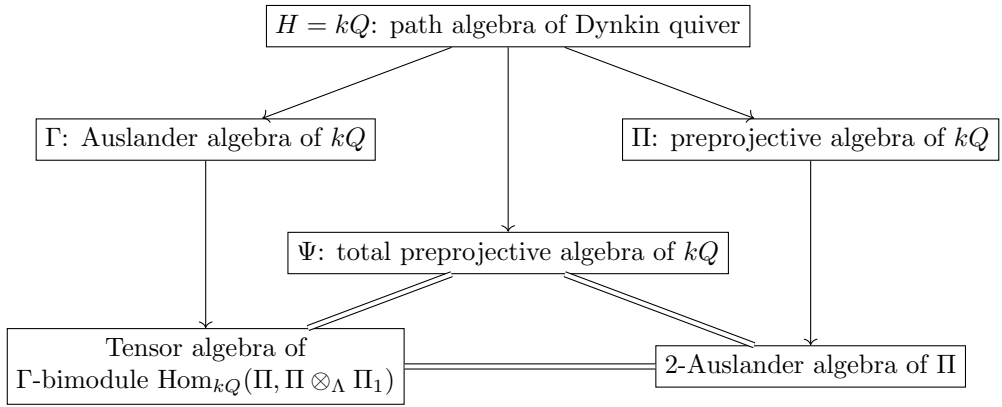
\[\begin{tikzpicture}[every text node part/.style={align=center}]
\tikzset{
    equ/.style={-,double equal sign distance},
  }
\node[draw] (vAusA) at (0,3.5) {$\Gamma$: Auslander algebra of $kQ$};
\node[draw] (vPsi) at (4,2) {$\Psi$: total preprojective algebra of $kQ$};
\node[draw] (vT) at (0,0.5) {Tensor algebra of \\$\Gamma$-bimodule $\Hom_{kQ}(\Pi,\Pi\otimes_\Lambda\Pi_1)$};
\node[draw] (vAusPi) at (8,0.5) {$2$-Auslander algebra of $\Pi$};
\node[draw] (vPi) at (8,3.5) {$\Pi$: preprojective algebra of $kQ$};
\node[draw] (v1) at (4,5) {$H=kQ$: path algebra of Dynkin quiver};
\draw  (v1) edge[->] (vPi);
\draw  (vPi) edge[->] (vAusPi);
\draw  (v1) edge[->] (vPsi);
\draw  (vPsi) edge[equ] (vAusPi);
\draw  (vAusPi) edge[equ] (vT);
\draw  (vPsi) edge[equ] (vT);
\draw  (vAusA) edge[->] (vT);
\draw  (v1) edge[->] (vAusA);
\end{tikzpicture}\]

Next, we give a quiver-with-relation presentation of the total preprojective algebras of $H=KQ$ for $Q$ of $\DynA\DynD\DynE$-type.
Let $\Gamma=\End_H(\Pi)$ be the Auslander algebra of $H$.  Suppose that we have a presentation of $\Gamma$ by quiver with relations $\Gamma \cong KQ_\Gamma/I_\Gamma$.
Then the vertices of $Q_\Gamma$ are in bijection to the indecomposable direct summands $X$ of $\Pi$.  
Define a new quiver $\widetilde{Q}_\Gamma$ by adding to the quiver $Q_\Gamma$ a new arrow $q_X \colon \tau^-(X) \rightarrow X$ for every indecomposable non-injective summand $X$ of $\Pi$.
Locally, we have the following subquiver in $\widetilde{Q}_{\Gamma}$.
\[\begin{tikzcd}
	X \dar[swap]{a} & {\tau^-(X)} \lar[swap]{q_X} \dar{\tau^-(a)} \\
	Y & {\tau^-(Y)} \lar{q_Y}
\end{tikzcd}\]
For each arrow $a \colon X \rightarrow Y$ of $Q_\Gamma$, let
\[r_a:=\left\{\begin{array}{ll}
q_X a-\tau^-(a) q_Y&\mbox{if $X$ and $Y$ are non-injective,}\\
q_X a&\mbox{if $X$ is non-injective and $Y$ is injective,}\\
0&\mbox{if $X$ is injective.}
\end{array}\right.\]

\begin{theorem}\label{mainpresentation} \cite[Theorem 1.7]{CIM2}
Let $H=KQ$ be the path algebra of simply-laced Dynkin type quiver $Q$ with $\Pi$ the direct sum of all indecomposable $H$-modules. Then we have a $K$-algebra isomorphism
\[\Psi\simeq  K\widetilde{Q}_\Gamma/(I_\Gamma,r_a\mid a\in (Q_\Gamma)_1).\]
\end{theorem}

We remark that this coincides with the combinatorial description given in \cite{GLS2} for such algebras.  We demonstrate the above recipe for quiver-and-relation in the following example.
\begin{example}
Let $H=KQ$ be the path algebra of the quiver
\[\begin{tikzcd}
	1 & 2 & 3
	\arrow[from=1-1, to=1-2]
	\arrow[from=1-2, to=1-3].
\end{tikzcd}\]
The Auslander algebra $\Gamma$ of $H$ is given by $\Gamma=KQ'/I'$ with
\[Q'=\begin{tikzcd}
	&& 3 \\
	& 2 && 5 \\
	1 && 4 && 6
	\arrow["d", from=1-3, to=2-4]
	\arrow["b", from=2-2, to=1-3]
	\arrow["c", from=2-2, to=3-3]
	\arrow["f", from=2-4, to=3-5]
	\arrow["a", from=3-1, to=2-2]
	\arrow["e", from=3-3, to=2-4]
\end{tikzcd}\;\text{ and }\; 
I'=\langle ac,ef,bd-ce\rangle.\]
Following the recipe, the total preprojective algebra $\Psi$ is isomorphic to $KQ''/I''$ with
\[Q'' = \begin{tikzcd}
	&& 3 \\
	& 2 && 5 \\
	1 && 4 && 6
	\arrow["d", from=1-3, to=2-4]
	\arrow["b", from=2-2, to=1-3]
	\arrow["c", from=2-2, to=3-3]
	\arrow["{q_2}", from=2-4, to=2-2]
	\arrow["f", from=2-4, to=3-5]
	\arrow["a", from=3-1, to=2-2]
	\arrow["e", from=3-3, to=2-4]
	\arrow["{q_1}", from=3-3, to=3-1]
	\arrow["{q_4}", from=3-5, to=3-3]
\end{tikzcd}\;\text{ and }\; I''=\langle ac,ef,bd-ce, q_1 a-e q_2,q_2 b,q_2c-fq_4,q_4 e\rangle.\]
This algebra has $K$-dimension 27 and is a higher Auslander algebra of global dimension $3$.
\end{example}

Apart from their natural description by quiver-and-relations, total preprojective algebras of Dynkin type enjoy several nice homological properties.
In a forthcoming work \cite{CCKMP}, we explore the Koszulity and quasi-heredity of total preprojective algebras, the description of the Koszul duals and their homological properties.

\section{Generalisations to higher dimensions}
In this last section, we generalise from hereditary algebras to a more general class of algebras.
Let $A$ be a finite dimensional algebra of global dimension $\leq d$ and set $\tau_d:=\tau \Omega^{d-1}$ the \emph{$d$-th higher Auslander-Reiten translate}. Since $A$ has global dimension $\leq d$, we have $\tau_d=D \Ext_A^d(-,A)$ and $\tau_d^{-1} = \Ext_A^d(DA,-)$.
An algebra $A$ with $\gldim A \leq d$ is said to be \emph{$\tau_d$-finite} if $\tau_d^i$ vanishes for $i \geq n$ for an integer  $n$, which is also equivalent to the vanishing of $\tau_d^{-i}$ for $i \geq n$. 
This notion generalises \emph{$d$-representation-finite algebras}, that is, an algebra $A$ that admits a $d$-cluster tilting module with $\gldim A \leq d$. 
We refer to \cite{I2,I3} for more information on $\tau_d$-finite and $d$-representation-finite algebras.
The next definition generalises (finite dimensional) preprojective algebras by replacing representation-finite hereditary algebras with $\tau_d$-finite algebras:
\begin{definition}
Let $A$ be a $\tau_d$-finite algebra.
Let $X \in \mod A$ be a direct summand of $\bigoplus\limits_{i\geq 0}^{}{\tau_d^{-i}(A)}$, then the \emph{higher preprojective algebra} of $X$ is defined as 
$$\Psi_A^X:= \bigoplus\limits_{i \geq 0}^{}{\Hom_A(X,\tau_d^{-i}(X))}.$$
When $X=\bigoplus\limits_{i\geq 0}^{}{\tau_d^{-i}(A)},$ we call $\Psi_A^X$, or for short just $\Psi$, the \emph{total preprojective algebra} of $A$ and when $X$ is projective, then we call $\Psi_A^X$ a \emph{higher contracted preprojective algebra.}
\end{definition}
The higher preprojective algebra of $A$ is the algebra $\Psi_A^A$, which we usually simply denote by $\Pi=\Pi(A)$.
Note that if $A$ is $d$-representation-finite, then as in the hereditary case we have $\Pi(A) \cong \bigoplus\limits_{i\geq 0}^{}{\tau_d^{-i}(A)}$ as a right $A$-module. For a description of higher preprojective algebras by quiver and relations in some special cases, we refer to \cite{GI}.

In general, $\Pi(A)$ is not selfinjective and the homological properties of $\Pi(A)$ are rather mysterious for general $\tau_d$-finite algebras. 
But when specialising to $d$-representation-finite algebras, classical results (somewhat) generalise:
\begin{theorem} \cite{IO}
Let $A$ be a $d$-representation-finite algebra with higher preprojective algebra $\Pi(A)$.
Then $\Pi(A)$ is a selfinjective and twisted periodic algebra.
\end{theorem}
Here, recall that a selfinjective algebra $B$ is called \emph{twisted periodic} if every simple $B$-module is periodic. It is an open problem, called the \emph{Periodicity Conjecture}, whether being twisted periodic implies that the algebra is periodic, we refer for example to \cite{CDIM} for more on this conjecture and a relation to another open problem on twisted fractionally Calabi-Yau algebras.

As a generalisation of Theorem \ref{explicit2CTobjectpreprojective}, we have:
\begin{theorem} \cite[Proposition 2.5]{CIM2}
Let $A$ be a $d$-representation-finite algebra with higher preprojective algebra $\Pi$.
Then $\Pi \otimes_A \Pi$ is a $(d+1)$-cluster-tilting object of $\mod \Pi$.
\end{theorem}

The next theorem gives the generalisation of Theorem \ref{main kQ} to $\tau_d$-finite algebras.
\begin{theorem}\label{main d}
Let $A$ be a finite dimensional algebra that is $\tau_d$-finite, $\Pi=\Pi(A)$ be the $(d+1)$-preprojective algebra of $A$, and $\Psi$ be the total preprojective algebra of $A$.
\begin{enumerate}[\rm(a)]
\item We have an isomorphism of algebras
\[\Psi\simeq\End_{\Pi}(\Pi\otimes_A\Pi).\]
\item If, moreover, $A$ is $d$-representation-finite, then $\Psi$ is a $(d+1)$-Auslander algebra, that is, $\gldim \Psi\le d+2\le \domdim \Psi$ holds.
\item Let $\Gamma:=\End_A(\Pi)$. Then $\Psi$ is isomorphic to the tensor algebra 
\[\Psi\simeq T_\Gamma\Hom_A(\Pi,\Pi\otimes_A\Pi_1).\]
\end{enumerate}
\end{theorem}


The recipe to calculate quiver-and-relations of total preprojective algebras of $d$-representation-finite algebras $A$ works exactly as in the $KQ$ case: First calculate the quiver-and-relations of the higher Auslander algebra $\Gamma=\End_A(\Pi)$. Then add arrows $q_X\colon \tau_d^{-1}(X) \rightarrow X$ for each indecomposable non-injective $X$ and add relation $r_a$ for each arrow $a$ of $\Gamma$.
We demonstrate this with an easy example:
\begin{example}
Let $A$ be the Nakayama algebra $KQ/I$ given by quiver and relations as follows:
$Q=$\[\begin{tikzcd}
	1 & 2 & 3 & 4
	\arrow["a", from=1-1, to=1-2]
	\arrow["b", from=1-2, to=1-3]
	\arrow["c", from=1-3, to=1-4]
\end{tikzcd}\]
$I= \langle abc \rangle$.  This is $2$-representation-finite where the $2$-cluster tilting module is $A\oplus DA$ with $\tau_2^{-1}(P_3)=I_1$ and $\tau_2^{-1}(P_4)=I_2$.
The higher preprojective algebra of $A$ is given by the quiver
\[\begin{tikzcd}
	1 & 2 & 3 & 4 \ar[lll, "t"', bend left]
	\arrow["a", from=1-1, to=1-2]
	\arrow["b", from=1-2, to=1-3]
	\arrow["c", from=1-3, to=1-4]
\end{tikzcd}\]
with generating relations $abc,bct,cta,tab$.
The higher Auslander algebra of $A$ is given by the quiver
\[\begin{tikzcd}
	1 & 2 & 3 & 4 & 5 & 6
	\arrow["a", from=1-1, to=1-2]
	\arrow["b", from=1-2, to=1-3]
	\arrow["c", from=1-3, to=1-4]
	\arrow["d", from=1-4, to=1-5]
	\arrow["e", from=1-5, to=1-6]
\end{tikzcd}\]
with generating relations $abc,bcd,cde$.
The total preprojective algebra of $A$ is given by the following quiver 
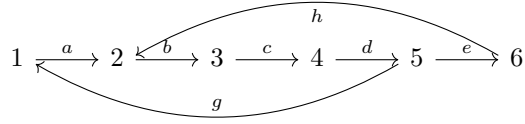
\[\begin{tikzcd}
	1 & 2 & 3 & 4 & 5 & 6
	\arrow["a", from=1-1, to=1-2]
	\arrow["b", from=1-2, to=1-3]
	\arrow["c", from=1-3, to=1-4]
	\arrow["d", from=1-4, to=1-5]
	\arrow["g"', shift right, bend left, from=1-5, to=1-1]
	\arrow["e", from=1-5, to=1-6]
	\arrow["h", shift left, bend right, from=1-6, to=1-2]
\end{tikzcd}\]
with generating relations $abc,bcd,cde,eh-ga,dg$.
\end{example}

For more complicated explicit examples of quiver-and-relations of higher preprojective algebras and higher total preprojective algebras, we refer to the example section in \cite{CIM2}.

We have now seen that the generalisations of our results from representation-finite path algebras to $d$-representation-finite algebras work smoothly for total preprojective algebras.
In contrast to that, we have not really explored the class of contracted preprojective algebras in the $d$-representation-finite case.
The study of higher contracted preprojective algebras seems to be much more complicated.
We pose the following problem:
\begin{problem}
For which $d$-representation-finite algebra $A$ is the following satisfied: the higher contracted preprojective algebra $\Psi_A^M$ is Cohen-Macaulay in the sense of Auslander and Reiten for all projective $A$-modules $M$?
\end{problem}
In Section 2, we saw that every path algebra $KQ$ with $Q$ Dynkin has the property that all contracted preprojective algebras are Cohen-Macaulay.  But the next example shows that, in general, higher contracted preprojective algebras are not Cohen-Macaulay, but one might find other nice homological properties exhibited by $\Psi_A^M$, such as being higher Auslander algebras.
\begin{example}
Consider the Auslander algebra of $KQ$, where $Q$ is of linear oriented Dynkin type $A_3$.
It is well-known that this is $2$-representation-finite as part of the higher Auslander algebras of type $A_n$ family, see \cite{I3}.
The higher preprojective algebra $A=KQ'/I'$ of the Auslander algebra of $KQ$ is given by quiver $Q'$ and admissible relations $I'$ as follows:

$Q'$=
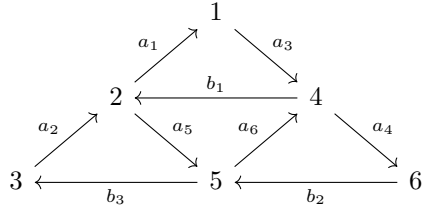
\[\begin{tikzcd}
	&& 1 \\
	& 2 && 4 \\
	3 && 5 && 6
	\arrow["{a_3}", from=1-3, to=2-4]
	\arrow["{a_1}", from=2-2, to=1-3]
	\arrow["{a_5}", from=2-2, to=3-3]
	\arrow["{b_1}"', from=2-4, to=2-2]
	\arrow["{a_4}", from=2-4, to=3-5]
	\arrow["{a_2}", from=3-1, to=2-2]
	\arrow["{a_6}", from=3-3, to=2-4]
	\arrow["{b_3}", from=3-3, to=3-1]
	\arrow["{b_2}", from=3-5, to=3-3]
\end{tikzcd}\]
$I'=\langle a_3 b_1,a_2 a_5,b_2 a_6,a_5 a_6-a_1 a_3,a_5b_3,a_4b_2-b_1a_5,b_1a_1,b_3a_2-a_6b_1,a_6a_4 \rangle.$
Now let $B_1:=eAe$ with $e=e_1+e_2+e_3+e_4+e_5$ be the higher contracted preprojective algebra. Then $B_1=KQ''/I''$ has quiver $Q''=$
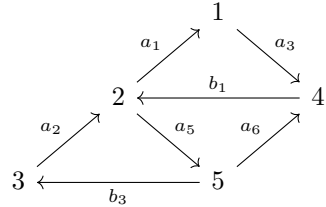
\[\quad \begin{tikzcd}
	&& 1 \\
	& 2 && 4 \\
	3 && 5
	\arrow["{a_3}", from=1-3, to=2-4]
	\arrow["{a_1}", from=2-2, to=1-3]
	\arrow["{a_5}", from=2-2, to=3-3]
	\arrow["{b_1}"', from=2-4, to=2-2]
	\arrow["{a_2}", from=3-1, to=2-2]
	\arrow["{a_6}", from=3-3, to=2-4]
	\arrow["{b_3}", from=3-3, to=3-1]
\end{tikzcd}\]
with relations $I''=\langle a_3b_1, a_2a_5, a_5 a_6-a_1 a_3, a_5b_3, b_3a_2-a_6b_1, b_1a_1 \rangle$.
$B_1$ is a higher Auslander algebra of global dimension 7.
Let $B_2:=f A f$ with $f$ the idempotent $f:=e_1+e_2+e_4+e_5$.
Then $B_2=KQ'''/I'''$ with quiver $Q'''=$
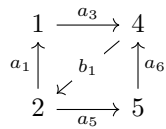
\[\begin{tikzcd}
	1 & 4 \\
	2 & 5
	\arrow["{a_3}", from=1-1, to=1-2]
	\arrow["{b_1}"{description}, from=1-2, to=2-1]
	\arrow["{a_1}", from=2-1, to=1-1]
	\arrow["{a_5}"', from=2-1, to=2-2]
	\arrow["{a_6}"', from=2-2, to=1-2]
\end{tikzcd}\]
and relations $I'''=\langle a_3b_1, a_1a_3-a_5a_6, b_1a_1, a_6b_1a_5 \rangle$.
$B_2$ has dominant dimension equal to three and is representation-finite with 14 indecomposable modules and has finitistic injective dimension and finitistic projective dimension equal to 3.
This implies, by \cite[Proposition 2.4]{CIM}, that the Ext-maximal module $W$ must be given by $W=\Omega^3(D B_2) \oplus U$, where $U$ is the direct sum of all indecomposable projective-injective $B_2$-modules.
Now $\dim B_2=12$, while $\dim \Hom_{B_2}(W,W)=14$. Thus $B_2$ is not Cohen-Macaulay as Cohen-Macaulay algebras $C$ with dualizing module $W$ satisfy that $\End_C(W) \cong C$, see \cite{AR2}.

\end{example}
\section{Further questions and problems}
We end this survey article with some further questions and problems.

We have studied higher preprojective algebras and total preprojective algebras for $d$-representation-finite algebras, but in general not much is known for the much larger class of $\tau_d$-finite algebras, which motivates our next problem:
\begin{problem}
Study higher preprojective algebras and total preprojective algebras for general $\tau_d$-finite algebras of global dimension $\leq d$. In particular, when are they Cohen-Macaulay, Iwanaga-Gorenstein or of finite global dimension?  
\end{problem}

One of the biggest open homological conjectures for finite dimensional algebras is the finitistic dimension conjecture, which predicts that the finitistic projective dimension is always finite.
We have mentioned that every Cohen-Macaulay algebra with dualizing module $W$ has finite finitistic dimension, and it is equal to $\idim W$ \cite{AR2}.
In particular, the finitistic dimension of contracted preprojective algebras of Dynkin type is at most 2 by Theorem \ref{CMmaintheorem}.
In fact one can show that all the algebras $\Psi_H^M$ for hereditary representation-finite algebras $H$ and arbitrary $M \in \mod H$ have finite finitistic projective dimension by using the following theorem due to Igusa and Todorov:
\begin{theorem} \cite[Corollary 8]{IT}
Let $A$ be a finite dimensional algebra of global dimension $\leq 3$. Then every algebra of the form $\End_A(P)$ has finite finitistic projective dimension for $P$ a projective $A$-module.
\end{theorem}

It follows that $\Psi_{KQ}^M$ has finite finitistic projective dimension for any $M\in \mod KQ$ using this result combined with our Theorem \ref{main kQ}, which states that total preprojective algebras of $\DynA\DynD\DynE$-type have global dimension $\leq 3$.

We do not know how to prove that general algebras of the form $\Psi_A^M$ have finite finitistic projective dimension and pose this as the next problem:
\begin{problem}
Let $A$ be a $\tau_d$-finite algebra of global dimension $\leq d$ and $M$ a direct summand of $\bigoplus\limits_{i\geq 0}^{}{\tau_d^{-i}(A)}$. Show that $\Psi_A^M$ has finite finitistic projective dimension.
\end{problem}

While we know that all algebras $\Psi_{KQ}^M$ have finite finitistic projective dimension for path algebras of Dynkin type $KQ$, it would be interesting to see how large the finitistic projective dimension can be for a given Dynkin type. 
For that reason we give the following definition here:
\begin{definition}
Let $H$ be a representation-finite hereditary algebra.
Define the \emph{total finitistic dimension} of $H$ as:
$$\operatorname{TFD}(H):= \sup \{ \fpdim \Psi_{H}^X \mid X \in \mod H \}.$$
\end{definition}

We pose the following problem:
\begin{problem}
Determine $\operatorname{TFD}(H)$ for representation-finite hereditary algebras $H$.
\end{problem}

Finally, let us show two interesting examples where $\operatorname{TFD}(KQ)>3$ and the preprojective algebra $\Psi_{KQ}^M$ for a $KQ$-module $M$ is in general not Cohen-Macaulay.  
These are computed using the GAP-package \cite{QPA}; for a brief guide (with preprojective algebra examples) on how it is used to calculate centralizer algebras of the form $eBe$ for a given quiver algebra $B$ and an idempotent $e$, we refer for example to \url{https://doi.org/10.5281/zenodo.17514267}. 
    \begin{example} \label{exoticexamples}
    Let $A=KQ$ with $Q$ given by the Dynkin quiver $Q=$
\[\begin{tikzcd}
	1 & 2 & 3 & 4
	\arrow[from=1-1, to=1-2]
	\arrow[from=1-2, to=1-3]
	\arrow[from=1-3, to=1-4]
\end{tikzcd}\]

The total preprojective algebra $\Psi$ is isomorphic to the algebra $KQ'/I'$ given by quiver and relations as follows:
$Q'=$
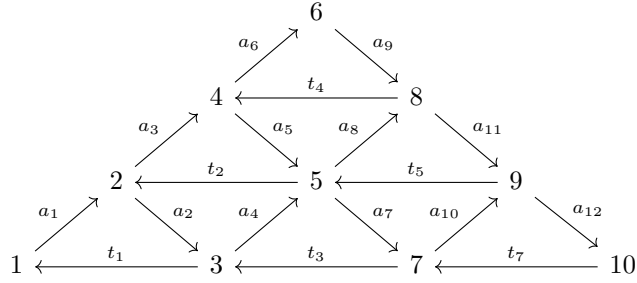
\[\begin{tikzcd}
	&&& 6 \\
	&& 4 && 8 \\
	& 2 && 5 && 9 \\
	1 && 3 && 7 && 10
	\arrow["{a_9}", from=1-4, to=2-5]
	\arrow["{a_6}", from=2-3, to=1-4]
	\arrow["{a_5}", from=2-3, to=3-4]
	\arrow["{t_4}"', from=2-5, to=2-3]
	\arrow["{a_{11}}", from=2-5, to=3-6]
	\arrow["{a_3}", from=3-2, to=2-3]
	\arrow["{a_2}", from=3-2, to=4-3]
	\arrow["{a_8}", from=3-4, to=2-5]
	\arrow["{t_2}"', from=3-4, to=3-2]
	\arrow["{a_7}", from=3-4, to=4-5]
	\arrow["{t_5}"', from=3-6, to=3-4]
	\arrow["{a_{12}}", from=3-6, to=4-7]
	\arrow["{a_1}", from=4-1, to=3-2]
	\arrow["{a_4}", from=4-3, to=3-4]
	\arrow["{t_1}"', from=4-3, to=4-1]
	\arrow["{a_{10}}", from=4-5, to=3-6]
	\arrow["{t_3}"', from=4-5, to=4-3]
	\arrow["{t_7}"', from=4-7, to=4-5]
\end{tikzcd}\]
$I'=\langle a_1 a_2, a_2 t_1, a_2 a_4- a_3 a_5, a_4 a_7, a_5 t_2, a_5 a_8-a_6 a_9, 
  a_7 a_{10}-a_8 a_{11}, a_9 t_4, a_{10} a_{12}, a_4 t_2-t_1 a_1, a_7 t_3-t_2 a_2, a_8 t_4-t_2 a_3, 
  a_{10}t_5-t_3 a_4, a_{11}t_5-t_4 a_5, a_{12} t_6- t_5 a_7 \rangle.$
Let $e$ be the idempotent $e=e_2+e_7+e_8+e_9$ and $B_1=e \Psi e$. Then $B_1 \cong KQ''/I''$
with $Q''$=
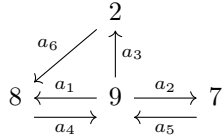
\[\begin{tikzcd}
	&& 2 \\
	{} & 8 & 9 & 7
	\arrow["{a_6}"', from=1-3, to=2-2]
	\arrow["{a_4}"', shift right=3, from=2-2, to=2-3]
	\arrow["{a_3}"', from=2-3, to=1-3]
	\arrow["{a_1}"', from=2-3, to=2-2]
	\arrow["{a_2}", from=2-3, to=2-4]
	\arrow["{a_5}", shift left=3, from=2-4, to=2-3]
\end{tikzcd}\]
and $I''=\langle a_1 a_4-a_2 a_5, a_4 a_3, a_5 a_2, a_6 a_4, a_3 a_6-a_1 a_4 a_1 \rangle.$

$B_1$ has global dimension equal to 6 and thus also finitistic projective dimension equal to 6. In fact we have TFD($KQ$)=6 for this choice of $Q$ as we have verified with the help of \cite{QPA}. Here we used that the delooping level (see \cite{Ge}) gives an upper bound for the finitistic injective dimension and the delooping level can be calculated for a given quiver algebra using \cite{QPA}.
Now choose $f=e_1+e_2+e_4+e_9$ and let $B_2:=f \Psi f$. Then $B_2 \cong K Q''' /I'''$ with 
$Q'''=$
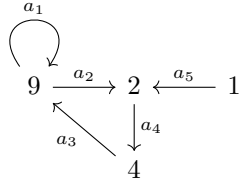
\[\begin{tikzcd}
	9 & 2 & 1 \\
	& 4
	\arrow["{a_1}", from=1-1, to=1-1, loop, in=55, out=125, distance=10mm]
	\arrow["{a_2}", from=1-1, to=1-2]
	\arrow["{a_4}", from=1-2, to=2-2]
	\arrow["{a_5}"', from=1-3, to=1-2]
	\arrow["{a_3}", from=2-2, to=1-1]
\end{tikzcd}\]
and relations $I'''= \langle a_1^2, a_1 a_2, a_3 a_1, a_3 a_2, a_4 a_3 \rangle.$
$B_2$ has 28 indecomposable modules and finitistic projective dimension equal to 1 and finitistic injective dimension equal to 2. Thus $B_2$ can not be Cohen-Macaulay, as Cohen-Macaulay algebras have finitistic projective dimension and finitistic injective dimension both equal to $\idim W$ for the dualizing module $W$, but finitistic projective and finitistic injective dimensions do not coincide.
\end{example}

\section*{Acknowledgements}
Many of the examples in this article have been found and verified with the help of the GAP-package \cite{QPA}. We thank the anonymous referee for helpful comments and corrections.

\end{document}